\documentclass[a4paper,11pt]{article}
\usepackage{graphicx}
\usepackage{a4wide}
\usepackage[english]{babel}
\usepackage{amsfonts}
\usepackage{amsmath}
\usepackage{amssymb}
\usepackage{dsfont}
\usepackage{url}
\usepackage{float}
\usepackage{amsthm,mathrsfs}
\newtheorem{lemma}{Lemma}
\newtheorem{theorem}{Theorem}

\newcommand {\E} {\mathbb{E}}

\DeclareMathOperator {\diag}{diag}

\newcommand {\p} {\mathbb{P}}

\newcommand {\R} {\mathbb{R}}

\makeatletter\def\blfootnote{\xdef\@thefnmark{}\@footnotetext}\makeatother
\allowdisplaybreaks
\title{\bf On a problem of Bourgain concerning the $L^1$-norm of exponential sums}
\author{Christoph Aistleitner\footnote{Department of Applied Mathematics, School of Mathematics and Statistics, University of New South Wales, Sydney NSW 2052, Australia. \mbox{e-mail}:
\texttt{aistleitner@math.tugraz.at}. The author is supported by an Erwin Schr\"odinger Fellowship of the Austrian Science Fund (FWF).}}
\begin{document}

\date{}
\maketitle

\blfootnote{{\bf MSC 2010:} 42A05, 33B10}
\blfootnote{{\bf keywords:} $L^1$ norm, exponential sums, trigonometric polynomials}

\begin{abstract}
Bourgain posed the problem of calculating 
$$
\Sigma = \sup_{n \geq 1} ~\sup_{k_1 < \dots < k_n} \frac{1}{\sqrt{n}} \left\| \sum_{j=1}^n e^{2 \pi i k_j \theta} \right\|_{L^1([0,1])}.
$$
It is clear that $\Sigma \leq 1$; beyond that, determining whether $\Sigma < 1$ or $\Sigma=1$ would have some interesting implications, for example concerning the problem whether all rank one transformations have singular maximal spectral type. In the present paper we prove $\Sigma \geq \sqrt{\pi}/2 \approx 0.886$, by this means improving a result of Karatsuba. For the proof we use a quantitative two-dimensional version of the central limit theorem for lacunary trigonometric series, which in its original form is due to Salem and Zygmund.
\end{abstract}

\section{Introduction and statement of results}

For any $n \geq 1$, set
$$
\Sigma_n = \sup_{k_1 < \dots < k_n} \frac{1}{\sqrt{n}} \left\| \sum_{j=1}^n e^{2 \pi i k_j \theta} \right\|_{1},
$$
where the supremum is taken over all sets of $n$ distinct positive integers, and the norm is the $L^1$ norm on the unit interval. Furthermore, set
$$
\Sigma = \sup_{n \geq 1} \Sigma_n.
$$
The problem of calculating $\Sigma_n$ and $\Sigma$ can be seen as the counterpart of the problem asking for optimal \emph{lower} bounds for the $L^1$ norm of trigonometric sums. Littlewood~\cite[p.~168]{hl} conjectured that there exists a constant $c$ such that for any distinct positive integers $k_1, \dots, k_n$ the inequality
\begin{equation} \label{li}
\left\| \sum_{j=1}^n e^{2 \pi i k_j \theta} \right\|_{1} \geq c \log n
\end{equation}
holds. This conjecture was confirmed in 1981 by Konyagin~\cite{kon} and, independently, McGehee, Pigno and Smith~\cite{mps}. It is easily seen that \eqref{li} is optimal (up to the value of the constant).\\

The problem of calculating $\Sigma$ appears in a paper of Bourgain~\cite{bourg}, in the context of his criterion for determining the maximal spectral type of rank one transformations. Deciding whether $\Sigma<1$ or $\Sigma=1$ would be a major step towards confirming or disproving the conjecture of Klemes and Reinhold~\cite{kr}, asserting that \emph{all} rank one transformations have singular spectral type. Bourgain~\cite[Proposition 2]{bourg} proved that
$$
\Sigma_n \leq 1 - \frac{c \log n}{n}.
$$
The strongest lower bound for $\Sigma$ until now is, as far as I know, $\Sigma \geq 1/\sqrt{2} \approx 0.707$, which is implicit in a Theorem of Karatsuba (in~\cite[Theorem 1]{karatsuba}, the minimal possible value of $K(N)$ is $2N(N-1)+N$, which implies $\sup_{N \to \infty} (\Delta(N))^{-1/2} N^{-1/2} = 1/\sqrt{2}$; this result is based on relating the $L^1$ norm of exponential sums to the number of solutions of certain Diophantine equations, see also~\cite{boch,emi,gar} for the context). The purpose of the present paper is to improve this lower bound for $\Sigma$. More precisely, our main result is the following theorem.

\begin{theorem} \label{th1}
$$
\Sigma \geq \frac{\sqrt{\pi}}{2} \approx 0.886.
$$
\end{theorem}

For the proof we use a quantitative two-dimensional version of the central limit theorem for lacunary trigonometric sums, which in its original form is due to Salem and Zygmund~\cite{sz}. More precisely, we will show that for the sequence $k_j = 8^j,~j \geq 1,$ we have
\begin{equation} \label{seq}
\lim_{n \to \infty} \frac{1}{\sqrt{n}} \left\| \sum_{j=1}^n e^{2 \pi i k_j \theta} \right\|_{1} = \frac{\sqrt{\pi}}{2}.
\end{equation}
Choosing our sequence $(k_j)_{j \geq 1}$ in such a way that it satisfies the Hadamard gap condition $k_{j+1}/k_j \geq q > 1$ for the growth factor $q = 8$ is just to simplify our proofs; any $q>1$ would be sufficient. Furthermore it is known that the central limit theorem for trigonometric series even holds for all sub-lacunary series satisfying
$$
\frac{k_{j+1}}{k_j} \geq 1 + \frac{c_j}{\sqrt{j}} \qquad \textrm{for $c_j \to \infty$}
$$
(Erd\H os~\cite{erdos}), and in a certain probabilistic sense it also holds for ``almost all'' slowly growing sequences (see Berkes~\cite{berkes} and Fukuyama~\cite{fuk1,fuk2} for details). This means that it is very likely that we could prove \eqref{seq} for many other sequences $(k_j)_{j \geq 1}$, and we take the specific sequence $k_j=8^j,~j \geq 1,$ just for the sake of simplicity. It should be noted that the central limit theorem for lacunary trigonometric series already appeared in the context of Bourgain's criterion for the spectral type of rank one transformations in a paper of El Abdalaoui~\cite{ela} (see also~\cite{aho}).\\

\section{Preliminaries}

\begin{lemma}[{\cite[Lemma 2]{steklov}}]   \label{lemmacf}
Let $P_1,P_2$ be probability measures on $\R^2$, and write $p_1,
p_2$ for the corresponding characteristic functions. Then for all $T_1,
T_2,\delta_1,\delta_2,x,y>0$
\begin{eqnarray*}
& & \left| P_1^*([-x,x] \times [-y,y]) - P_2^*([-x,x] \times [-y,y]) \right| \\
& \leq & xy \int_{(s,t) \in
[-T_1,T_1] \times [-T_2,T_2]} |p_1(s,t)-p_2(s,t)| ~d(s,t) \\
& & + xy  \left(  \delta_1^{-1} \delta_2~ \exp\left(- T_1^2
\delta_1^2/2 \right) + \delta_1 \delta_2^{-1} ~\exp \left(-
T_2^2 \delta_2^2/2 \right) \right),
\end{eqnarray*}
where
$$
P_1^* = P_1 \star H, \quad P_2^* = P_2 \star H,
$$
and $H$ is a two-dimensional normal distribution with density
$$
(2 \pi \delta_1 \delta_2)^{-1} ~e^{-\delta_1^2 u^2/2 -
\delta_2^2 v^2/2}.
$$
\end{lemma}

We will use the classical trigonometric formulas
\begin{eqnarray} 
\sin x \sin y & = & \frac{\cos(x-y) - \cos(x+y)}{2}, \nonumber\\
\cos x \cos y & = & \frac{\cos(x-y)+\cos (x+y)}{2}, \label{trig}\\
\sin x \cos y & = & \frac{\sin(x-y)+\sin(x+y)}{2}.\nonumber
\end{eqnarray}
We will also use the fact that for a random variable $Z$ having two-dimensional normal distribution with expectation vector $(0,0)$ and covariance matrix $\diag(\sigma^2,\sigma^2)$ we have
\begin{equation} \label{moment}
\E |Z| = \sqrt{\frac{\pi  \sigma^2}{2}}.
\end{equation}

\section{Proof of Theorem~\ref{th1}}

Our proof is based on the method of Salem and Zygmund in~\cite{sz}, utilizing characteristic functions, but additionally considering sine- and cosine-functions simultaneously and finally, using the smoothing inequality in Lemma~\ref{lemmacf}, giving a quantitative estimate for the deviation between the distribution of our normalized lacunary sum and the normal distribution. It would also be possible to prove the convergence of the distribution of the normalized sums $\sum e^{2 \pi i k_j \theta}$ to a complex normal distribution, but since working with complex distributions is rather unusual and would not yield any simplification or shortening, we prefer to use the two-dimensional real argument instead.\\
For one-dimensional lacunary series, quantitative versions of the Salem-Zygmund central limit theorem have already been proved, based on an approximation by martingale differences and an almost sure invariance principle of Strassen (see~\cite{abe} for details). Recently this martingale approach has been extended to a multi-dimensional setting by Moore and Zhang~\cite{moore} for proving the law of the iterated logarithm for lacunary series. However, for our purpose the classical method of using characteristic functions seems to be much easier.\\

Let $n \geq 1$ be given. We set
$$
k_j = 8^j, \qquad 1 \leq j \leq n,
$$
and, for $\theta \in [0,1]$,
$$
\mu(\theta) = \frac{\sum_{j=1}^n \sin 2 \pi k_j \theta}{\sqrt{n}}, \qquad \nu(\theta) = \frac{\sum_{j=1}^n \cos 2 \pi k_j \theta}{\sqrt{n}}.
$$
In the following paragraph, we will understand the symbol ``$\E$'' with respect to the probability space $([0,1]),\mathcal{B}([0,1]),\lambda)$, where $\lambda$ is the Lebesgue measure. We write $\exp(x)$ for $e^x$. For $s,t \in \R$ set
$$
\varphi (s,t) = \E (e^{i s \mu + i t \nu}).
$$
Using the classical equality
$$
e^{ix} = (1+ix) e^{-x^2/2 + w(x)}, \qquad \textrm{where $|w(x)| \leq x^3$},
$$
we get
\begin{eqnarray*}
& & e^{i s \mu + i t \nu} \\
& = & \left(\prod_{j=1}^n \exp\left(\frac{is \sin 2 \pi k_j \theta}{\sqrt{n}} \right)\right) \cdot \left( \prod_{j=1}^n \exp\left( \frac{it\cos 2 \pi k_j \theta}{\sqrt{n}} \right) \right) \\
& = &  \underbrace{\left(\prod_{j=1}^n \left(1 + \frac{i s \sin 2 \pi k_j \theta}{\sqrt{n}} \right) \right) \cdot \left(\prod_{j=1}^n \left(1 + \frac{it \cos 2 \pi k_j \theta}{\sqrt{n}} \right) \right)}_{=: \alpha(s,t)} \times \nonumber\\
& & \quad \times \exp \left( \sum_{j=1}^n \left( -\frac{(s \sin 2 \pi k_j \theta)^2 + (t \cos 2 \pi k_j \theta)^2}{2n} + w\left(\frac{s \sin 2 \pi k_j \theta}{\sqrt{n}} \right) +  w\left(\frac{t \cos 2 \pi k_j \theta}{\sqrt{n}} \right) \right) \right) \\
& = & \alpha(s,t) ~\exp \left( - \frac{s^2+t^2}{4} + \underbrace{\sum_{j=1}^n \left(\frac{(s^2-t^2)\cos 4 \pi k_j \theta}{4n} + w\left(\frac{s \sin 2 \pi k_j \theta}{\sqrt{n}} \right) +  w\left(\frac{t \cos 2 \pi k_j \theta}{\sqrt{n}} \right) \right)}_{=: \beta(s,t)} \right) \\
& = & \alpha(s,t) \exp \left( -\frac{s^2+t^2}{4} + \beta(s,t) \right).
\end{eqnarray*}
Note that
\begin{eqnarray}
|\alpha(s,t)| & = & \left(\prod_{j=1}^n \left( 1 + \frac{(s \sin 2 \pi k_j \theta)^2}{n} \right)^{1/2} \right) \cdot \left(\prod_{j=1}^n \left( 1 + \frac{(t \cos 2 \pi k_j \theta)^2}{n} \right)^{1/2} \right) \nonumber\\
& \leq & \prod_{j=1}^n \exp \left( \frac{(s \sin 2 \pi k_j \theta)^2 + (t \cos 2 \pi k_j \theta)^2}{2n} \right) \nonumber\\
& = & \exp \left( \frac{s^2+t^2}{4} + \sum_{j=1}^n \frac{(t^2-s^2) \cos 4 \pi k_j \theta}{4n}\right). \label{alp}
\end{eqnarray}
Furthermore, by construction of the sequence $(k_j)_{1 \leq j \leq n}$ we have
\begin{equation} \label{ga}
\E \left( \alpha(s,t) \right) = 1.
\end{equation}
To see this, note that we have
\begin{equation*}
\alpha(s,t) = \sum_{\substack{(\delta_1, \dots, \delta_n) \in \{0,1\}^n,\\(\hat{\delta}_1, \dots, \hat{\delta}_n)\in\{0,1\}^n}} \left( \left( \prod_{j=1}^n  \left( \frac{is \sin 2 \pi k_j \theta}{\sqrt{n}} \right)^{\delta_j} \right) \cdot \left( \prod_{j=1}^n \left( \frac{it \cos 2\pi k_j \theta}{\sqrt{n}} \right)^{\hat{\delta}_j} \right) \right), \label{cho}
\end{equation*}
and using the equations in \eqref{trig} to transform a product
\begin{equation} \label{sumr}
\prod_{j=1}^n \left(\left( is \sin 2 \pi k_j \theta \right)^{\delta_j} \left( it \cos 2\pi k_j \theta\right)^{\hat{\delta}_j}\right)
\end{equation}
into a sum of trigonometric functions it is easily seen that by the construction of the sequence $(k_1, \dots, k_n)_{1 \leq j \leq n}$ only sine-functions and cosine-functions with nonzero frequency appear in this sum representation of \eqref{sumr}, which means that
$$
\E \left( \prod_{j=1}^n \left( \left( is \sin 2 \pi k_j \theta \right)^{\delta_j} \left( it \cos 2\pi k_j \theta\right)^{\hat{\delta}_j} \right) \right) = 0
$$
unless $(\delta_1, \dots, \delta_n) = (\hat{\delta}_1,\dots, \hat{\delta}_n) = (0,\dots,0)$.\\

Thus 
\begin{eqnarray}
& & \left| \varphi(s,t) - e^{-(s^2+t^2)/4} \right| \nonumber\\
& = & \left| \E \left(\alpha(s,t) ~\exp \left(- \frac{s^2+t^2}{4} + \beta(s,t) \right) -  e^{-(s^2+t^2)/4} \right) \right| \label{g1}\\
& = & \left| \E \left(\alpha(s,t) \left(\exp \left(- \frac{s^2+t^2}{4} + \beta(s,t) \right) -  e^{-(s^2+t^2)/4} \right) \right) \right|  \label{g2} \\
& = & \left| \E \left(\alpha(s,t) e^{-(s^2+t^2)/4} \left(e^{\beta(s,t)} - 1\right) \right) \right| \nonumber\\
& \leq & \E \left(\left|\alpha(s,t)\right|~e^{-(s^2+t^2)/4} ~\left| e^{\beta(s,t)} - 1 \right| \right) \label{alpa}\\
& \leq & \E \left( \left| \exp\left( \sum_{j=1}^n \left( w\left(\frac{s \sin 2 \pi k_j \theta}{\sqrt{n}} \right) +  w\left(\frac{t \cos 2 \pi k_j \theta}{\sqrt{n}} \right) \right)\right)  - \exp \left( \sum_{j=1}^n \frac{(t^2-s^2) \cos 4 \pi k_j \theta}{4n}\right) \right| \right) \nonumber\\
& \leq & \exp \left( \frac{|s|^3+|t|^3}{n^{1/2}} \right) - 1 + \E \left( \left| \exp \left( \sum_{j=1}^n \frac{(t^2-s^2) \cos 4 \pi k_j \theta}{4n} \right) - 1 \right| \right). \label{deviation}
\end{eqnarray}
Here \eqref{ga} was used to get from \eqref{g1} to \eqref{g2}, and \eqref{alp} was used to estimate \eqref{alpa}. It is easily seen that by construction of $(k_j)_{1 \leq j \leq n}$ we have
$$
\E \left( \left( \sum_{j=1}^n \cos 4 \pi k_j \theta \right)^4 \right) \leq n^2,
$$
which by Markov's inequality gives 
$$
\lambda \left( x \in [0,1]:~\left| \sum_{j=1}^n \cos 4 \pi k_j x \right| \geq n^{3/4} \right) \leq 1/n,
$$
and consequently
\begin{eqnarray} 
& & \E \left( \left| \exp \left( \sum_{j=1}^n \frac{(t^2-s^2) \cos 4 \pi k_j \theta}{4n} \right) - 1 \right| \right) \nonumber\\
& \leq & \exp \left( n^{-1/4} (s^2+t^2) \right) - 1 + \exp \left(s^2+t^2\right) n^{-1}.\label{deviation2}
\end{eqnarray}

Now let $X$ be a two-dimensional random vector having the same distribution as $(\mu(\theta),\nu(\theta)),~\theta \in [0,1]$, let $Y$ be a random vector having two-dimensional normal distribution with expectation vector $(0,0)$ and covariance matrix $\diag(1/2,1/2)$, and let $Z$ be a random vector having two-dimensional normal distribution with expectation $(0,0)$ and covariance matrix $\diag((\log n)^{-1/8},(\log n)^{-1/8})$, such that $Z$ is independent of $X$ and $Y$. Using \eqref{deviation} and \eqref{deviation2} to estimate the deviation between $\varphi(s,t)$ and the characteristic function of $Y$, and using Lemma~\ref{lemmacf} for $T_1=T_2=(\log n)^{1/4}$, we obtain that for any $x,y \in \R, ~|x|\leq (\log n)^{1/4}, ~|y| \leq (\log n)^{1/4}$,
\begin{eqnarray*}
& & \Big| \p \left( X+Z \in [-x,x] \times [-y,y]\right) - \p \left( Y+Z \in [-x,x] \times [-y,y] \right) \Big| \\
& \leq & (\log n) \left( \exp \left( \frac{2 (\log n)^{3/8}}{n^{1/2}} \right) - 1 + \exp \left(\frac{2 (\log n)^{1/2}}{n^{1/4}} \right) - 1 + \frac{\exp \left( 2 (\log n)^{1/2} \right)}{n} \right) \\
& & + 2 (\log n)^{1/2} e^{-(\log n)^{1/4}/2} \\
& \leq & c_1 (\log n)^{-1}
\end{eqnarray*}
for some appropriate constant $c_1$. This implies
$$
\E \left(|Y+Z| \cdot \mathds{1}_{|Y+Z| \leq (\log n)^{1/4}}\right) \leq \E \left(|X+Z| \cdot \mathds{1}_{|X+Z| \leq (\log n)^{1/4}}\right) + c_1 (\log n)^{-3/4}.
$$
Consequently, using \eqref{moment} and the fact that $Y+Z$ is also normally distributed with expectation $(0,0)$ and covariance $\diag(1/2+(\log n)^{-1/8},1/2+(\log n)^{-1/8})$, we finally obtain
\begin{eqnarray*}
\frac{1}{\sqrt{n}} \left\| \sum_{j=1}^n e^{2 \pi i k_j \theta} \right\|_1 & = & \left\| \mu(\theta) + i \nu(\theta) \right\|_1 \\
& = & \E |X| \\
& \geq & \E |X+Z| - \E|Z| \\
& \geq & \E \left(|X+Z| \cdot \mathds{1}_{|X+Z| \leq (\log n)^{1/4}}\right) - \E|Z| \\
& \geq & \E \left(|Y+Z| \cdot \mathds{1}_{|Y+Z| \leq (\log n)^{1/4}}\right) - c_1 (\log n)^{-3/4} - \E |Z| \\
& \geq & \E |Y+Z| - \E \left(|Y+Z| \cdot \mathds{1}_{|Y+Z| > (\log n)^{1/4}} \right) - c_1 (\log n)^{-3/4} - \E |Z| \\
& \geq & \frac{\sqrt{\pi}}{2} - c_2 (\log n)^{-1/16}.
\end{eqnarray*}
for some appropriate positive constant $c_2$. Thus we have established \eqref{seq}, which proves Theorem~\ref{th1}.

\end{document}